\documentclass{amsart}
\usepackage{amsfonts}
\usepackage{amsmath}
\usepackage{amssymb}
\usepackage[utf8]{inputenc}

\newtheorem{theorem}{Theorem}
\newtheorem{corollary}{Corollary}
\newtheorem{lemma}{Lemma}
\newtheorem{proposition}{Proposition}

\newtheorem{definition}{Definition}
\newtheorem{remark}{Remark}

\begin{document}

\title[Identities and central polynomials of real graded division algebras]{Identities and central polynomials of real graded division algebras}

\author{Diogo Diniz}
\address{Unidade Acadêmica de Matemática, Universidade Federal de Campina Grande, Campina Grande, PB, 58429-970, Brazil} \thanks{The first author acknowleges support by CNPq grant \# 303534/2013-3}
\email{diogo@mat.ufcg.edu.br}

\author{Claudemir Fidelis}
\address{Unidade Acadêmica de Matemática, Universidade Federal de Campina Grande, Campina Grande, PB, 58429-970, Brazil}
\email{diogo@mat.ufcg.edu.br}

\author{Sérgio Mota}\address{ Universidade Estadual de Santa Cruz - UESC\\
Rodovia Ilhéus/Itabuna, Km 16\\ 45662-900 Ilhéus - BA, Brazil} \email{smalves@uesc.br}

\begin{abstract}
Let $A$ be a finite dimensional real algebra with a division grading by a finite abelian group $G$. In this paper we provide finite basis for the $T_G$-ideal of graded identities and for the $T_G$-space of graded central polynomials for $A$.
\end{abstract}
\maketitle

\section{Introduction}

Let $K$ be a field and $A=\oplus_{g\in G} A_g$ a $K$-algebra graded by the group $G$, the set of graded polynomial identities for $A$ is an ideal of the free $G$-graded algebra $K\langle X_G \rangle$ invariant under graded endomorphisms, such ideals are called $T_G$-ideals. A central problem in the theory of algebras satisfying polynomial identities (p. i. algebras) is the Specht problem on the existence of a finite basis for the ideal of identities for any algebra in a given class (Associative, Lie, Jordan, etc.). This problem was solved by A.Kemer \cite{K} for (ungraded) associative algebras over a field of characteristic zero. The corresponding result for p.i. algebras graded by a finite group $G$ was proved by I. Sviridova \cite{S} (for abelian groups) and by E. Aljadeff and A. Kanel-Belov \cite{AK}. Explicit finite basis are known in a few cases only, for the (ungraded) identities of the matrix algebra of order $2$ a finite basis was determined by Yu. P. Razmyslov \cite{Ra}, a minimal basis was latter obtained by V. Drensky \cite{Dr}.
A polynomial in $K\langle X_G \rangle$ is a graded central polynomial for $A$ if the result of every admissible substitution is a central element. The set of central polynomials for $A$ is a subspace of $K\langle X_G \rangle$ invariant under graded endomorphisms, such subspaces are called $T_G$-subspaces. Central polynomials play an important role in p.i. theory, we refer the reader to \cite{DrF}, \cite{Y} for more details and applications.

In his study of projective planes M. Hall proved that a division algebra for which a certain polynomial (the Hall polynomial, see Theorem \ref{Dr}) is an identity is either a field or a generalized quaternion algebra over its center. Latter I. Kaplansky generalized this result and proved that the conclusion holds if the division algebra satisfies any polynomial identity. Division gradings on simple real finite-dimensional algebras were classified by Yu. Bahturin and M. Zaicev \cite{BZ} and A. Rodrigo-Escudero \cite{R}. More generally finite-dimensional graded division algebras over the field of real numbers were recently classified in \cite{BZ2}. In this paper we provide finite basis for the $T_G$-ideal of graded identities and the $T_G$-space of graded central polynomials for the graded division algebras in  \cite{BZ}, \cite{R} and \cite{BZ2}.

If an algebra has a division grading, with support $G$, such that the component associated to the neutral element is one dimensional then every component is one dimensional. In this case homogeneous elements commute up to a scalar and the $G$-grading is regular (see Definition \ref{def.rel}). The division algebras in \cite{BZ}, \cite{R} and \cite{BZ2} are either regular, the tensor product of an algebra with a regular grading  and a basic algebra (presented in sections \ref{divgrad} and \ref{basic2}) or a complex division grading on a matrix algebra over $\mathbb{C}$ viewed as a real algebra (Pauli grading). In \cite{DM} a basis for graded identities of the tensor product of a graded algebra $A$ and an algebra with  a regular grading is obtained from a basis for the graded identities for $A$. In Section \ref{section.AtensorR} we prove the analogous result for graded central polynomials, we remark that a particular case of this result was proved in \cite{ABK}. Therefore a basis for the graded identities and central polynomials of the tensor product division algebras is obtained once a basis for the basic algebras is known. For most of the basic algebras such basis are known, in sections \ref{TPDG} and \ref{basic2} we obtain the description for the remaining basic algebras. In Section \ref{NRPG} we consider the Pauli gradings that are not regular. 

\section{Preliminaries}

Let $K$ be a field of characteristic zero, $G$ a group with identity element $e$. A \textit{grading by the group $G$} (or a $G$-grading) on the vector space $V$ is a decomposition $\Gamma: V=\oplus_{g\in G} V_g$. Let $\Gamma^{\prime}:W=\oplus_{h\in H} W_h$ be an $H$-graded vector space, a morphism of graded vector spaces is a linear map $f:V\rightarrow W$ such that for every $g\in G$ there exists $h\in H$ such that $f(V_g)\subseteq W_h$. If $f$ has an inverse which is also a morphism of graded vector spaces we say that the gradings $\Gamma$ and $\Gamma^{\prime}$ are \textit{equivalent}. A $G$-grading on an algebra $A$ is a grading on $A$ as a vector space, $\Gamma: A=\oplus_{g\in G} A_g$, such that $A_{g}A_{h}\subseteq A_{gh}$ for every $g,h$ in $G$. The subspaces $A_{g}$ are called the \textit{homogeneous components} of $A$. An element $a\in \cup_{g\in G} A_g$ is called a \textit{homogeneous element}, if $a$ is non-zero there exists a unique $g\in G$ such that $a\in A_g$, this element is called the \textit{degree of $a$} (in the $G$-grading) and denoted $\mathrm{deg_G}\ (a)$. The \textit{support} of $\Gamma$ is the set $\mathrm{supp}\ \Gamma:=\{g\in G \hspace{0.1cm} | \hspace{0.1cm} A_g\neq 0\}$. Let $\Gamma^{\prime}: A=\oplus_{h\in H} A_h$ be another grading by $H$ on $A$, we say that $\Gamma^{\prime}$ is a \textit{coarsening} of $\Gamma$ if for every $g\in G$ there exists $h\in H$ such that $A_g\subseteq A_h$. Let $\Gamma:A=\oplus_{g\in G}$ and $\Gamma^{\prime}:B=\oplus_{h\in H}B_h$ be algebras graded by the groups $G$ and $H$ respectively. The gradings $\Gamma$ and $\Gamma^{\prime}$ are \textit{equivalent} if there exists an isomorphism of algebras $\varphi:A\rightarrow B$ that is an equivalence of graded vector spaces. The graded algebras $A$ and $B$ are \textit{weakly isomorphic} if there exists a group isomorphism $\alpha:G\rightarrow H$ and an algebra isomorphism $\varphi:A \rightarrow B$ such that $\varphi_{A_g}=B_{\alpha(g)}$ for every $g$ in $G$. If $G$ and $H$ are the supports $\mathrm{supp}\ A$ and $\mathrm{supp}\ B$, respectively, and if the grading $\Gamma$ is strong (i.e. $A_gA_h=A_{gh}$ for every $g,h \in G$) then the gradings $\Gamma$ and $\Gamma^{\prime}$ are equivalent if and only if they are weakly isomorphic (see \cite[Section 2]{BZ2}). The tensor product $A\otimes B$ has a canonical $G\times H$ grading where $(A\otimes B)_{(g,h)}=A_g\otimes B_h$.  The (graded) algebra $\Gamma: A=\oplus_{g\in G} A_g$ is a \textit{graded division algebra} (or $\Gamma$ is a division grading on $A$) if $A$ has a unit and every non-zero homogeneous element is invertible.

We denote by $K\langle X_G \rangle$ the free $G$-graded algebra, freely generated by the set $X_{G}=\{x_{i,g}|i\in \mathbb{N}, g\in G\}$. This algebra has a natural grading by $G$ where, for any $g$ in $G$, the homogeneous component $(K\langle X_G \rangle)_g$ is the subspace generated by the monomials $x_{i_1,g_1}\cdots x_{i_m,g_m}$ such that $g_{1}\cdots g_m=g$. Let $f(x_{1,g_1},\dots, x_{m,g_m})$ be a polynomial in $K\langle X_G \rangle$, an $m$-tuple $(a_1,\dots, a_m)$ such that $a_i\in A_{g_i}$ for $i=1,\dots,m$ is called an \textit{$f$-admissible substitution} (or simply an admissible substitution). If $f(a_1,\dots,a_m)=0$ for every admissible substitution $(a_1,\dots, a_m)$ we say that the polynomial $f$ is a \textit{graded polynomial identity} for $A$.  We denote by $Z(A)$ the centre of $A$ and we say that $f$ is a \textit{graded central polynomial} for $A$ if $f(a_1,\dots, a_m)\in Z(A)$ for every admissible substitution $(a_1,\dots, a_m)$. We denote by $Id_G(A)$ and $C_G(A)$ the set of all all graded identities and all graded central polynomials for $A$, respectively. We say that the polynomial $f$ is an \textit{proper graded central polynomial} for $A$ if $f\in C_G(A)\setminus Id_G(A)$, in other words $f$ is $G$-graded central polynomial but not a $G$-graded identity for $A$.

It is clear that $Id_G(A)$ and $C_G(A)$ are subspaces of $K\langle X_G \rangle$ invariant under all endomorphisms of $K\langle X_G \rangle$, such subspaces are called \textit{$T_G$-spaces}. Moreover, the set $Id_G(A)$ of all graded identities for $A$ is an ideal of $K\langle X_G \rangle$ invariant under all endomorphisms, such ideals are called \textit{$T_G$-ideals}. The intersection of $T_G$-ideals (respect. $T_G$-spaces) of $K\langle X_G \rangle$ is also a $T_G$-ideal (respect. $T_G$-space). A subset $S\subset K\langle X_G \rangle$ \textit{generates the $T_G$-ideal $Id_G(A)$ (resp. $T_G$-space $C_G(A)$)} if it equals the intersection of all $T_G$-ideals (resp. $T_G$-spaces) in $K\langle X_G \rangle$ that contain $S$. If $G=\{e\}$ we recover the definition of ordinary polynomial identities and central polynomials, in this case we omit the sub index $G$ and the sub indices that are elements of $G$ and  use the notation $K\langle X \rangle$ for the free associative algebra and $x_{i}$ for the variables. We denote $[x_1,x_2]=x_1x_2-x_2x_1$ the commutator of $x_1$ and $x_2$. The \textit{standard polynomial} of degree $n$ is the polynomial $$S_n(x_1,\dots, x_n)=\sum_{\sigma \in S_n}(-1)^{\sigma} x_{\sigma(1)}\cdots x_{\sigma(n)},$$ where $(-1)^{\sigma}$ denotes the sign of the permutation $\sigma$.

Let $A=\oplus_{g\in G}$ and $B=\oplus_{h\in H}B_h$ be weakly isomorphic algebras graded by the groups $G$ and $H$ respectively and let $\alpha:G\rightarrow H$ and $\varphi:A \rightarrow B$ be a group isomorphism and an algebra isomorphism, respectively, such that $\varphi_{A_g}=B_{\alpha(g)}$ for every $g$ in $G$. A polynomial $f\in K\langle X_G \rangle$ is a graded identity for $A$ (resp. graded central polynomial) if and only if $\Phi (f)$ is a graded identity (resp. graded central polynomial) for $B$, where $\Phi: K\langle X_G \rangle \rightarrow K\langle X_H \rangle$ is the isomorphism such that $\Phi (x_{ig})=x_{i\varphi(g)}$. Moreover if $S$ is a basis for the $T_G$-ideal $Id_G(A)$ (resp. $T_G$-space $C_G(A)$) then $\Phi(S)$ is a basis for the $T_G$-ideal $Id_G(B)$ (resp. $T_G$-space $C_G(B)$).

Let $\Gamma: A=\oplus_{g\in G} A_g$ be a $G$-graded division algebra, the support $H=\mathrm{supp}\ \Gamma$ is a subgroup of $G$ and $A=\oplus_{h\in H} A_h$ is an $H$-grading on $A$. If $S\subset K\langle X_H \rangle$ is a basis for $Id_H(A)$ (resp. $C_H(A)$) then the set $S\cup \{x_{1g} \hspace{0.1cm | \hspace{0.1cm}} g \in G\setminus H\}$ is a basis for $Id_G(A)$ (resp. $C_G(A)$).

\begin{remark}
Henceforth whenever $\Gamma:A=\oplus_{g\in G} A_g$ is a $G$-graded division algebra we assume that $G=\mathrm{supp}\ \Gamma$.
\end{remark}

Next we recall the definition of a regular grading (see \cite{AO}, \cite{RS}).
\begin{definition}\label{def.rel}
Let $R$ be an algebra over the field $K$ and \[
R=\bigoplus_{h\in H} R_h,
\]
a grading by the abelian group $H$. The $H$-grading above on $R$ is said to be regular if there exists a commutation function $\beta:H\times H \rightarrow K^{\times}$ such that
\begin{enumerate}
  \item [(P1)] For every $n$ and every $n$-tuple $(h_1,\dots, h_n)$ of elements of $H$, there exist $r_1,\dots,r_n$ such that $r_j\in R_{h_j}$, and $r_1\cdots r_n\neq0$.
  \item [(P2)] For every $h_1,h_2 \in H$ and for every $a_{h_1}\in A_{h_1}$, $a_{h_2}\in A_{h_2}$, we have
      \[
      a_{h_1}a_{h_2}=\beta(h_1,h_2)a_{h_2}a_{h_1}.
      \]
\end{enumerate}
The regular $H$-grading on $R$ is said to be minimal if for any $h\neq e$ there exists $h^{\prime}$ such that $\beta(h,h^{\prime})\neq 1$.
\end{definition}

In the definition above the commutation function $\beta:H\times H \rightarrow K^{\times}$ is a \textit{skew-symmetric bicharacter}, i.e., for every $h_0 \in H$ the maps $h\mapsto \beta(h_0,h)$ and $h\mapsto \beta(h,h_0)$ are characters of the group $H$ and $\beta(h_2,h_1)=\beta(h_1,h_2)^{-1}$ for every $h_1,h_2 \in H$. We say that the bicharacter $\beta$ is \textit{non-degenerate} if $\beta(h_0,h)=1$ for every $h\in H$ implies that $h_0=e$.

\section{Graded Identities and Central Polynomials: From $A$ To $A\otimes R$}\label{section.AtensorR}

Let $R$ be an algebra with a regular grading by the abelian group $H$ and $A$ an algebra graded by the group $G$. Our main result in this section, Theorem \ref{main_Theom}, provides a basis for the central polynomials of the algebra $A\otimes R$ with its canonical $G\times H$-grading obtained from a basis for the $G$-graded central polynomials of $A$ that contains a basis for $Id_G(A)$ as a $T_G$-space.

Let $\textbf{g}=(g_1,\cdots, g_n)$ be an $n$-tuple of elements of the group $G$. We denote by $P_{\bf g}$ the set of polynomials multilinear in the variables $x_{1,g_1},\cdots, x_{n,g_n}$. We fix sequence ${\bf h}=(h_1,\dots,h_n)$ of elements of $H$. The description of the basis for $A\otimes R$ is in terms of multilinear maps $\phi_{\textbf{h}}: P_{\bf g}\rightarrow P_{\bf g\times h}$, where $\textbf{g}\times \textbf{h}=((g_1,h_1),\cdots, (g_n,h_n))$, introduced in \cite{BDr}. A permutation $\tau\in S_n$, together with $\textbf{h}$, determines a non-zero scalar $\lambda_{\tau}^{\textbf{h}}$ such that
\begin{equation}\label{lambda}
r_1\cdots r_n = \lambda_{\tau}^{\textbf{h}} r_{\tau(1)}\cdots r_{\tau(n)},
\end{equation}
whenever $r_i \in R_{h_i}$ for $i=1,\dots,n$. We define $$\phi_{\textbf{h}}(x_{g_{\tau(1)},\tau(1)}\cdots x_{g_{\tau(n)},\tau(n)})=\lambda_{\tau}^{\textbf{h}}x_{(g_{\tau(1)},h_{\tau(1)}),\tau(1)}\cdots x_{(g_{\tau(n)},h_{\tau(n)}),\tau(n)},$$ and extend to $P_{\bf g}$ by linearity.

\begin{theorem}\cite[Theorem 5.5]{DM}\label{basis} 
Let $R$ be an algebra with a regular grading by the abelian group $H$ and let $A$ be any $G$-graded algebra. If $S\subseteq \bigcup_{\substack{{\bf g}\in G^n \\ n\in \mathbb{N}} } P_{\bf g}$ is a basis, of multilinear polynomials, for the $T_G$-ideal $Id_G(A)$, then the set
\[\widetilde{S}=\{\phi_{\bf h}(f)\,\hspace{0,1cm}|\hspace{0,1cm}\, f(x_{1g_1},\dots, x_{ng_n})\in S,\, {\bf h}\in H^n\},\] is a basis for the $T_{G\times H}$-ideal $Id_{G\times H}(A\otimes R)$.
\end{theorem}

Let $\beta:H\times H \rightarrow K^{\times}$ denote the bicharacter associated to the regular $H$-grading on $R$. Denote $$H^{\prime}=\{h^{\prime}\in H \hspace{0,1cm}|\hspace{0,1cm} \beta(h^{\prime},h)=1 \hspace{0,1cm} \forall h \in H\},$$ it is clear that $R_{h^{\prime}}\subseteq Z(R)$ for every $h^{\prime}\in H^{\prime}$. In this section we assume that $R$ has minimal center in the sense that
\begin{equation}\label{center}
Z(R)=\oplus_{h^{\prime}\in H^{\prime}}R_{h^{\prime}}.
\end{equation}

\begin{remark}\label{equal}
Let $f(x_{g_1,1},\cdots, x_{g_n,n})$ be a multilinear polynomial and $\textbf{h}$ an $n$-tuple in $H^{n}$. It follows from the definition of the map $\phi_{\bf h}$ (see \cite[Theorem 3.1]{BDr}) that  the equality
\begin{equation}\label{f}
\phi_{\bf h}(f)(a_1\otimes r_1,\cdots, a_n\otimes r_n)=f(a_1,\dots, a_n)\otimes r_1\cdots r_n,
\end{equation}
holds for every $\phi_{\bf h}(f)$-admissible substitution $(a_1\otimes r_1, \cdots, a_n\otimes r_n)$. This implies that $\phi_{\bf h}(f)$ is a graded identity for $A\otimes R$ if and only if $f$ is a graded identity for $A$. If equality (\ref{center}) holds then $\phi_{\bf h}(f)$ is a proper central polynomial for $A\otimes R$ if and only if $f$ is a proper central polynomial for $A$ and $(h_1\cdots h_n) \in H^{\prime}$.
\end{remark}

Our next theorem is the analogous result for graded central polynomials. In its proof we need the following lemma.

\begin{lemma}\cite[Lemma 5.4]{DM}\label{10}
Let $f=f(x_{g_11},\dots,x_{g_mm})\in P_{\bf g}$, and  $w_1,\dots,w_m$ monomials in $K\langle X_G\rangle$ such that $\mathrm{deg_G}\ (w_i)=g_i$ for each $i$ and $f(w_1,\dots,w_m)\in P_{\bf k}$, for some $\textbf{k} \in G^{n}$. If ${\bf h}\in H^n$, there exist ${\bf h'}\in H^m$ and homogeneous elements $b_1,\dots, b_m\in K\langle X|G\times H\rangle$ such that
$\phi_{\bf h}(f(w_1,\dots,w_m))=\gamma
\phi_{\bf h'}(f)(b_1,\dots,b_m)$.
\end{lemma}

\begin{theorem}\label{main_Theom}
Let $A$ be an algebra with a grading by a group $G$ and let $R$ be an algebra with a regular grading by an abelian group $H$, with bicharacter $\beta:H\times H \rightarrow K^{\times}$, such that $Z(R)=\bigoplus_{h\in H^\prime}R_h$, where $H^{\prime}=\{h^{\prime}\in H\hspace{0,1cm}| \hspace{0,1cm}\beta(h^{\prime},h)=1 \forall h \in H\}$. If $S=S_1\cup S_2 \subseteq \bigcup_{\substack{{\bf g}\in G^n \\ n\in \mathbb{N}} } P_{\bf g}$ is a basis, of multilinear polynomials, for the $T_G$-space $C_{G}(A)$ such that $S_1$ generates $Id_G(A)$ as a $T_G$-space then the set $\widetilde{S}=\widetilde{S_1}\cup \widetilde{S_2}$, where $$\widetilde{S_1}=\{\phi_{\bf h}(f)\hspace{0,1cm}|\hspace{0,1cm}f(x_{g_1,1},\dots, x_{g_n,n})\in S_1, \textbf{h}\in H^{n}\}$$ and $$\widetilde{S_2}=\{\phi_{\bf h}(f)\hspace{0,1cm}|\hspace{0,1cm}f(x_{g_1,1},\dots, x_{g_n,n})\in S_2, \textbf{h}=(h_1,\dots, h_n)\in H^{n}, h_1\cdots h_n\in H^{\prime}\}$$ is a basis for the $T_{G\times H}$-space $C_{G\times H}(A\otimes R)$.
\end{theorem}
\textit{Proof.}
The inclusion $\widetilde{S}\subseteq C_{G\times H}(A\otimes R)$ follows from Remark \ref{equal}, therefore $\langle \widetilde{S} \rangle_{T_{G\times H}} \subseteq C_{G\times H}(A\otimes R)$, where $\langle \widetilde{S} \rangle_{T_{G\times H}}$ denotes the $T_{G\times H}$-space generated by $S$. Let $f(x_{(g_1,h_1),1},\dots, x_{(g_n,h_n),n})$ be a multilinear element of $C_{G\times H}(A\otimes R)$. We denote $\textbf{g}=(g_1,\dots, g_n)$ and $\textbf{k}=((g_1, h_1),\dots, (g_n,h_n))$. It is clear from the definition that the map $\phi_{\bf h}$ is bijective. This implies that there exists a multilinear polynomial $f^{\prime}(x_{g_1,1},\dots, x_{g_n,n})$ in $P_{\bf g}$ such that $f=\phi_{\bf h}(f^{\prime})$.  It follows from Remark \ref{equal} that $f^{\prime}$ is a graded central polynomial for $A$. Therefore there exists polynomials $p_1,\dots, p_n$ in $S$,  $p_i$-admissible substitutions $(w_1^i,\dots,w_{k_i}^i)$ and scalars $\alpha_i$ such that
\begin{equation}\label{11}
f^{\prime}=\sum_{i=1}^n\alpha_ip_i(w_1^{i},\dots, w_{k_i}^{i}).
\end{equation}
Since $f^{\prime}$ and the polynomials $p_i$ are multilinear we may assume without loss of generality that the $w_k^{i}$ are multilinear monomials. Lemma \ref{10} implies that  there exists ${\bf h^i}\in H^{k_i}$, $\phi_{\bf h^i}(p_i)$-admissible substitutions $(b_1^i,\dots,b_{k_i}^i)$ and scalars $\gamma_i\in K^{\times}$ such that
\[\phi_{\bf h}(p_i(w_1^i,\dots,w_{k_i}^i))=\gamma_i\phi_{\bf h^i}(p_i)(b_1^i,\dots,b_{k_i}^i),\]
for $i=1,\dots, n$. It follows from (\ref{11}) and the above equalities that
\begin{equation}\label{ee}
f=\phi_{\bf h}(f^{\prime})=\sum_i{\gamma}^{\prime}_i\phi_{\bf h^i}(p_i)(b_1^i,\dots,b_{k_i}^i),
\end{equation}
where ${\gamma}^{\prime}_i=\alpha_i{\gamma}_i$. If $f$ is a graded identity for $A\otimes R$ then it follows from Remark \ref{equal} that $f^{\prime}$ is a graded identity for $A$. In this case we may assume that the polynomials $p_1,\dots, p_n$ lie in $S_1$, therefore the polynomials $\phi_{\bf h^i}(p_i)$ lies in $\widetilde{S_1}$ (which is a subset of $\widetilde{S}$) and equality (\ref{ee}) implies that $f$ lies in $\langle \widetilde{S} \rangle_{T_{G\times H}}$. Now we assume that $f$ is not a graded identity for $A\otimes R$. Since $f$ is multilinear there exists an admissible substitution of the form $(a_1\otimes r_1,\dots, a_n\otimes r_n)$ such that $a_1,\dots, a_n$ and $r_1,\dots, r_n$ are homogeneous elements in $A$ and $R$ respectively and $f(a_1\otimes r_1,\dots, a_n\otimes r_n)\neq 0$. Therefore $f(a_1\otimes r_1,\dots, a_n\otimes r_n)$ is a non-zero central element. Equality (\ref{f}) implies that $$f(a_1\otimes r_1,\dots, a_n\otimes r_n)=f^{\prime}(a_1,\dots, a_n)\otimes r_1\cdots r_n,$$ therefore $r_1\cdots r_n$ is a homogeneous non-zero central element in $R$ of degree $h_1\cdots h_n$. Since $R$ satisfies the equality $Z(R)=\bigoplus_{h\in H^\prime}R_h$ we conclude that $h_1\cdots h_n\in H^{\prime}$. The polynomials $p_1,\dots, p_n$ lie in $S$ and since $h_1\cdots h_n\in H^{\prime}$ we conclude that $\phi_{\bf h^i}(g_i)$ lies in $\widetilde{S}$ for $i=1,\dots, n$. Therefore equality (\ref{ee}) implies that $f$ lies in $\langle \widetilde{S} \rangle_{T_{G\times H}}$.
\hfill $\Box$

Note that $S=S_1\cup S_2$, where $S_1=\{x_1[x_2,x_3]x_4\}$ and $S_2=\{x_1\}$, is a basis of multilinear polynomials for the central polynomials of the field $K$ such that $S_1$ generates the identities of $K$ as a $T$-space. Let $R$ be an algebra with a regular $H$-grading satisfying the hypothesis of the previous theorem. The algebra $K\otimes R$ with its $\{e\}\times H$-grading is weakly isomorphic to $R$ with its $H$-grading. Hence the previous theorem yields the following:

\begin{corollary}\label{descr.alg.reg}
Let $R$ denote an algebra with a regular grading by the abelian group $H$, with corresponding bicharacter $\beta_R:H\times H \rightarrow K^{\times}$, such that $Z(R)=\bigoplus_{h\in H^\prime}R_h$, where $H^{\prime}=\{h^{\prime}\in H\hspace{0,1cm}| \hspace{0,1cm}\beta(h^{\prime},h)=1 \forall h \in H\}$. The polynomials
\begin{equation}\label{central.reg}
 x_{1 h}, h\in H^\prime
\end{equation}
and
\begin{equation} \label{id.reg}
x_{1 h_1}\left(x_{2 h_2}x_{3 h_3}-\beta_R(h_2,h_3)x_{3 h_3}x_{2 h_2}\right)x_{h_4,4}, \mbox{ where } h_i\in H, i=1,\dots, 4.
\end{equation}
form a basis for the $T_H$-space $C_{H}(R)$.
\end{corollary}

\begin{remark}\label{basisrem}
Analogously,  as a corollary of Theorem \ref{basis}, we conclude that the set $$x_{1,h_1}x_{2,h_2}-\beta(h_1,h_2)x_{2,h_2}x_{1,h_1},$$ where $h_1,h_2 \in H$ is basis for the $H$-graded identities of $R$.
\end{remark}

\begin{remark}
Note that $S_2=\{x_1\}$ is a basis for the $T$-space of central polynomials of the  field $K$, however no set of indeterminates $x_{ih}$ generates $T_H$-space $C_H(R)$. Hence as noted in \cite{ABK} the relation between a basis of $C_G(A)$ and $C_{G\times H}(A\otimes R)$ is not direct.
\end{remark}

\section{Basis for the Graded Identities and Central Polynomials for Division Gradings on Simple Real Algebras}\label{section.Application}

Finite dimensional simple real algebras with a division grading by a finite abelian group $G$ were classified in \cite{BZ} and \cite{BD}. In this section we present the main result of \cite{BD} and describe basis for the $T_G$-ideal of graded identities and the $T_G$-space of graded central polynomials of these algebras.

\subsection{Division Gradings on Simple Real Algebras}\label{divgrad}

Next we present the basic division gradings that are the building blocks of this classification. One such building block is the $\mathbb{Z}_2$-grading on the field $\mathbb{C}$ of complex numbers: $$\Gamma:\mathbb{C}=\mathbb{C}_0\oplus \mathbb{C}_1,$$ where $\mathbb{C}_0=\langle 1 \rangle$ and $\mathbb{C}_1=\langle i \rangle$. We will denote $\mathbb{C}^{(2)}$ the above grading on $\mathbb{C}$. The other building blocks are division gradings on the algebras $M_2(\mathbb{R})$, $\mathbb{H}$ and $M_2(\mathbb{C})$.

\subsubsection{Division gradings on $M_2(\mathbb{R})$}
Let $G=(a)_2\times(b)_2\cong\mathbb{Z}_2\times\mathbb{Z}_2$. A division grading on $S=M_2(\mathbb{R})$ may be obtained by means of Sylvester matrices below:
$$A=\left(
             \begin{array}{cc}
               1 & 0 \\
               0 & -1\\
             \end{array}
           \right)
, \hspace{0.1cm} B=\left(
             \begin{array}{cc}
               0 & 1 \\
               1 & 0 \\
             \end{array}
           \right), \hspace{0.1cm}C=\left(
             \begin{array}{cc}
               0 & 1 \\
               -1 & 0 \\
             \end{array}
           \right).$$
We have $A^2=B^2=I,AB=C=-BA$, therefore a division grading by $G$ is obtained where the homogeneous components are
$$S_{e}=\langle I\rangle,\hspace{0.1cm} S_a=\langle A\rangle,\hspace{0.1cm} S_b=\langle B\rangle, \hspace{0.1cm}S_c=\langle C\rangle,$$
here $c=ab$ and $I$ is the identity matrix. This grading is regular and the center of $M_2(\mathbb{R})$ is the neutral component. We denote this graded division algebra by ${M}_2^{(4)}$. The decomposition $$R_{e}=\langle I, C\rangle, \hspace{0.1cm} R_a=\langle A, B,\rangle$$ where $a$ is the generator of the group $H=(a)_2\cong \mathbb{Z}_2$, is a coarsening of $M_2^{(4)}$ which is also a division grading on $R=M_2(\mathbb{R})$. This grading will be denoted by ${M}_{2}^{(2)}$.

\subsubsection{Division gradings on $\mathbb{H}$}

Let $\{1,i,j,k\}$ be the canonical basis of the quaternion algebra $\mathbb{H}$. The multiplication is determined by the relations $$i^2=j^2=-1, \hspace{0.1cm}ij=-ji=k.$$ Since $\mathbb{H}$ is a division algebra, any grading makes it a graded division algebra. There is, up to weak isomorphism, two non trivial gradings on $\mathbb{H}$ that we describe next. A grading by the group  $G=(a)_2\times (b)_2$ is the decomposition
$$R_{e}=\langle 1\rangle, \hspace{0.1cm}R_a=\langle i\rangle, \hspace{0.1cm}R_b=\langle j\rangle, \hspace{0.1cm} R_c=\langle k\rangle$$
where $c=ab$. This grading is regular and the neutral component is the center of the quaternions. We denote this graded division algebra by $\mathbb{H}^{(4)}$. A coarsening of this grading is
$$S_{e}=\langle 1,i\rangle, S_a=\langle j, k\rangle,$$ this is a grading by the group $H=(a)_2\cong \mathbb{Z}_2$ on $S=\mathbb{H}$. This grading will be denoted by $\mathbb{H}^{(2)}$.

\subsubsection{Division gradings on $M_2(\mathbb{C})$}

Let $G=(a)_4\times (b)_2$. A division $G$-grading on $R=M_2(\mathbb{C})$, which we denote $M_2^{(8)}$, is obtained with homogeneous components
\begin{eqnarray*}
&&R_{e}=\langle I \rangle, \hspace{0.1cm} R_{a}=\langle \omega A\rangle,\hspace{0.1cm} R_{a^{2}}=\langle iI \rangle,\hspace{0.1cm} R_{a^3}=\langle i\omega A \rangle\\
&&R_{b}=\langle C \rangle, \hspace{0.1cm}R_{ab}=\langle \omega B\rangle, \hspace{0.1cm}R_{a^{2}b}=\langle iC \rangle, \hspace{0.1cm}R_{a^3b}=\langle i\omega B \rangle,
\end{eqnarray*}
where $w=\frac{1}{\sqrt{2}}(1+i)$ is the 8-th root of 1 such that $\omega^2=i$.
A division grading on $S=M_2(\mathbb{C})$ by the group $H = (a)_4$ obtained as a coarsening of $M_2^{(8)}$ is
$$S_{e}=\langle I,C\rangle, S_a=\langle \omega A,\omega B\rangle, S_{a^2}=\langle iI,iC\rangle, S_{a^3}=\langle i\omega A,i\omega B\rangle,$$
this grading will be denoted by $M_2(\mathbb{C},\mathbb{Z}_4)$.

\subsubsection{Pauli gradings}
Now we describe the \textit{Pauli gradings} on the real algebra $R=M_n(\mathbb{C})$. The algebra $R$ is also an algebra over $\mathbb{C}$. If the decomposition $\Gamma: R=\oplus_{g\in G} R_g$ is a division grading by the group $G$ for $R$ as an algebra over $\mathbb{C}$ then it is also a division grading for $R$ over $\mathbb{R}$. Such gradings are called Pauli gradings. The division grading $\Gamma$ is a regular grading for $R$ as an algebra over $\mathbb{C}$ and can be described in terms of the associated bicharacter $\beta:G\times G \rightarrow \mathbb{C}^{\times}$. This grading is regular for $R$ as a real algebra if and only if the image of $\beta$ lies in $\mathbb{R}^{\times}$. We refer to $\Gamma$ as a \textit{non-regular Pauli grading} if it is not a regular grading for $R$ as a real algebra. We remark that $\Gamma$ is a regular Pauli grading if and only if $G\cong \mathbb{Z}_2^{2k}$ for some $K$.

Our study of the graded identities of graded division algebras relies on the classification of division gradings on finite dimensional simple real algebras, this classification was done in \cite{BZ} and \cite{R}. The main result of \cite{BZ} is the following:

\begin{theorem}\label{theor.BZ}
Any division grading on a real simple algebra $M_n(D)$, $D$ a real division algebra, is weakly isomorphic to one of the following types:
\begin{enumerate}
  \item[$D=\mathbb{R}$:]  \begin{enumerate}
                                \item [(i)] $({M}_{2}^{(4)})^{\otimes k}\cong C^{gr}(k,k)$;
                                \item [(ii)] ${M}_{2}^{(2)}\otimes({M}_{2}^{(4)})^{\otimes(k-1)}$, a coarsening of {\rm (i)};
                                \item [(iii)]${M}_{4}^{(4)}\otimes({M}_{2}^{(4)})^{\otimes(k-2)}$, a coarsening of {\rm (i)};
                              \end{enumerate}
															\vspace{0.3cm}
  \item[$D=\mathbb{H}$:]  \begin{enumerate}
                                \item [(iv)] $\mathbb{H}^{(4)}\otimes ({M}_{2}^{(4)})^{\otimes k}\cong C^{gr}(k+1,k-1)$;
                                \item [(v)] ${\mathbb{H}}^{(2)}\otimes({M}_{2}^{(4)})^{\otimes k}$, a coarsening of {\rm (iv)};
                                \item [(vi)] $\mathbb{H}\otimes({M}_{2}^{(4)})^{\otimes k}$, a coarsening of {\rm (v)};
                              \end{enumerate}
															\vspace{0.3cm}
  \item[$D=\mathbb{C}$:]  \begin{enumerate}
                                \item [(vii)] $\mathbb{C}^{(2)}\otimes ({M}_{2}^{(4)})^{\otimes k}\cong C^{gr}(k+1,k)$;
                                \item [(viii)] $\mathbb{C}^{(2)}\otimes{M}_{2}^{(2)}\otimes ({M}_{2}^{(4)})^{\otimes(k-1)}$, a coarsening of {\rm (vii)};
                                \item [(ix)] $\mathbb{C}^{(2)}\otimes\mathbb{H}\otimes ({M}_{2}^{(4)})^{\otimes k}$, a coarsening of {\rm (vii)};
                                \item [(x)] $ {M}_{2}^{(8)}\otimes({M}_{2}^{(4)})^{\otimes(k-1)}$;
                                \item [(xi)] $M_2(\mathbb{C},\mathbb{Z}_{4})\otimes({M}_{2}^{(4)})^{\otimes(k-1)}$, a coarsening of {\rm (x)};
                                \item [(xii)] ${M}_{2}^{(8)}{\otimes}\mathbb{H}\otimes ({M}_{2}^{(4)})^{\otimes(k-2)}$, a coarsening of {\rm (x)};
                                \item [(xiii)] Pauli gradings.
                              \end{enumerate}
\end{enumerate}
None of the gradings of different types or of the same type but with different values of $k$ is weakly isomorphic to the other.
\end{theorem}

\subsection{Graded Identities and Central Polynomials for Division Gradings on Simple Real Algebras}\label{basissec}

In this section we describe basis for the graded identities and central polynomials for the algebras in Theorem \ref{theor.BZ}. The gradings are either regular, a tensor product grading of one of the  building blocks presented in Section \ref{divgrad} and an algebra with a regular grading or a non-regular Pauli grading.

\subsubsection{Regular Division Gradings}\label{reg}

The algebras in Theorem \ref{theor.BZ} with a regular division grading are $(i)$, $(iv)$, $(vii)$, $(x)$ and the regular Pauli gradings in $(xiii)$. The tensor product of two algebras with regular gradings has a regular grading, we have the following remark.

\begin{remark}\label{tensorregular}
Let $A,B$ be algebras with regular gradings by the abelian groups $G$ and $H$, respectively. The tensor product grading on $A\otimes B$ is regular and the corresponding bicharacter $\beta_{A\otimes B}:(G\times H)\times (G\times H) \rightarrow K^{\times}$ is given by $$\beta_{A\otimes B}((g_1,g_2),(h_1,h_2))=\beta_A(g_1,g_2)\beta_B(h_1,h_2),$$ for $g_1,g_2\in G$ and $h_1h_2\in H$.
\end{remark}

If $\beta$ is the bicharacter corresponding to a regular grading on $R$ by the finite abelian group $G$, generated by the finite set $S=\{s_1,\dots, s_k\}$, then it is determined by the set of values $\beta(s_i,s_j)$, for $1\leq i<j \leq k$. Moreover each $\beta(s_i,s_j)$ is a root of the unit and therefore equals $\pm 1$. Hence $\beta$ is determined by the set of pairs $(s_i,s_j)$ such that $\beta(s_i,s_j)=-1$. One determines this pairs for the building blocks in Section \ref{divgrad} by inspection, we then use Remark \ref{tensorregular} to obtain the bicharacters for the algebras $(i)$, $(iv)$, $(vii)$, $(x)$ in Theorem \ref{theor.BZ}. The algebras in $(i)$ and $(iv)$ are canonically graded by the group $G_{2k}=(a_1)_2\times (b_1)_2\times \cdots \times (a_k)_2\times (b_k)_2$, and the algebras in $(vii)$ and $(x)$ by the groups $G_{2k+1}=(a_0)_2\times(a_1)_2\times (b_1)_2\times \cdots \times (a_k)_2\times (b_k)_2$, and $H_{k-1}=(a_1)_4\times (b_1)_2\times \cdots \times (a_k)_2\times (b_k)_2$ respectively. The bicharacters are presented in the following table:

\vspace{0.3 cm}

\begin{center}\begin{tabular}{|c|c|c|}
  \hline
  \textbf{Algebra} & \textbf{Grading group} & \textbf{Bicharacter}\\ \hline
  $(M_2^{(4)})^{\otimes k}$, $\mathbb{H}^{(4)}\times (M_2^{(4)})^{\otimes (k-1)}$ & $G_{2k}$ & $\beta_k(a_i,b_i)=-1$, $i=1,\dots, k$\\ \hline
  $\mathbb{C}^{(2)}\otimes (M_2^{(4)})^{\otimes k}$ & $G_{2k+2}$  & $\gamma_k(a_i,b_i)=-1$, $i=1,\dots, k$\\ \hline
	$M_2^{(8)}\times (M_2^{(4)})^{\otimes (k-1)}$ & $H_{2k-1}$ & $\delta_k(a_i,b_i)=-1$, $i=1,\dots, k$\\
  \hline
\end{tabular}
\end{center}

\vspace{0.3 cm}

A Pauli grading on $M_n(\mathbb{C})$ by the group $G$ is regular if and only if $G\cong \mathbb{Z}_2^{2k}$ and $n=2^{k}$ for some $k$. In this case the grading is weakly isomorphic to $M_2(\mathbb{C})^{\otimes k}$ where $M_2(\mathbb{C})$ has its canonical Pauli grading. This algebra can be canonically endowed with a grading by the group $(a_1)_2\times (b_1)_2\times \cdots \times (a_k)_2\times (b_k)_2$ in such a way that the corresponding bicharacter is $\beta_k$.

Note that the bicharacters $\beta_k$ and $\delta_k$ are non-degenarate and $$H^{\prime}=\{h^{\prime}\in H\hspace{0,1cm}| \hspace{0,1cm}\gamma_k(h^{\prime},h)=1 \forall h \in H\}=(a_0)_2.$$ Therefore a basis for the graded central polynomials of finite dimensional simple real algebras with a regular division grading are obtained from Corollary \ref{descr.alg.reg}.

\subsubsection{Tensor Product Division Gradings}\label{TPDG}
A basis for the graded identities and central polynomials for the algebras in $(ii)$, $(iii)$, $(v)$, $(vi)$, $(viii)$, $(ix)$, $(xi)$ and $(xii)$, of Theorem \ref{theor.BZ},  is obtained from Theorem \ref{basis} and Theorem \ref{main_Theom} once a suitable basis for the central polinomials of the algebras $\mathbb{H}$, $\mathbb{H}^{(2)}, {M}^{(2)}_{2}$ and $M_2(\mathbb{C}, \mathbb{Z}_4)$ is known. Such a basis is known for the algebras $\mathbb{H}$, $\mathbb{H}^{(2)}$ and ${M}^{(2)}_{2}$. 

\begin{theorem}\cite{Dr}\label{Dr}
Let $K$ be a field of characteristic zero. The standard polynomial $S_4(x_1, x_2, x_3, x_4)$ and the Hall polynomial $[[x_1, x_2]^2, x_3]$ form a basis for the identities of $M_2(K)$.
\end{theorem}

\begin{theorem}\cite{O}
Let $K$ be an infinite field of characteristic zero. The polynomials $[x_1,x_2][x_3,x_4]+[x_3,x_4][x_1,x_2]$ and $x_5S_4(x_1, x_2, x_3, x_4)$ form a basis for the central polynomials of $M_2(K)$.
\end{theorem}

\begin{remark}\label{hall}
It is well known that the algebras $\mathbb{H}$ and $M_2(\mathbb{R})$ satisfy the same polynomial identities, hence the previous theorems provide basis for the polynomial identities and central polynomials for $\mathbb{H}$.
\end{remark}

It was proved in \cite{BD} that $\mathbb{H}^{(2)}$ and ${M}^{(2)}_{2}$ satisfy the same graded polynomial identities as $R=M_2(\mathbb{R})$ with the $(a)_2$-grading 
\begin{equation}\label{elem}
R_{e}=\left(
\begin{array}{cc}
	\mathbb{R} & 0\\
	0 & \mathbb{R}
\end{array}\right) \hspace{0,1cm} R_{a}=\left(
\begin{array}{cc}
		0 & \mathbb{R}\\
		\mathbb{R} & 0
\end{array}\right).
\end{equation}
Basis for the graded identities and central polynomials for $R$ were determined in \cite{DV} and in \cite{BP}. 

\begin{lemma}\cite[Lemma 2]{DV} \label{2}
The set consisting of the polynomials $x_{1 e}x_{2 e}-x_{2 e}x_{1 e}$ and $x_{1,a}x_{2,a}x_{3,a}-x_{3,a}x_{2,a}x_{1,a}$ is a basis for the $T_{G^{\prime}}$-ideal of graded identities of $R$ with the grading in (\ref{elem}).
\end{lemma}

\begin{theorem}\cite[Theorem 4]{BP} \label{3}
The set of polynomials \[x_{1a}^2, \hspace{0,3cm} x_{1g}[x_{2,e},  x_{3e}]x_{4g}, \hspace{0,3cm} x_{1g}\left(x_{2a}x_{3a}x_{4a}-x_{4a}x_{3a}x_{2b}\right)x_{5g},\] where $g\in (a)_2$ generates the $T_{(a)_2}$-space $C_{(a)_2}(R)$.
\end{theorem}

\begin{remark}
It is well known that the polynomial $x_{1a}^2$ in the theorem above may be replaced by the multilinear polynomial $x_{1a}x_{2a}+x_{2a}x_{1a}$ to obtain a basis of multilinear polynomials.
\end{remark}

To complete our description we determine a basis for the graded identities of the algebra $M_2(\mathbb{C}, \mathbb{Z}_4)$. The basis is obtained using the following proposition.

\begin{proposition}\label{1}
Let $A=\oplus_{g\in G} A_g$ be an algebra graded by the group $G$ and assume that there exists an element $g$ of the group $G$ such that $A_g$ contains an invertible central element. If $Q:K\langle X_G \rangle \rightarrow K\langle X_{G/\langle g \rangle}\rangle$ is the homomorphism such that $Q(x_{i,g})=x_{i,q(g)}$, where $q:G\rightarrow G/\langle g \rangle$ is the canonical homomorphism, and $S^{\prime}$ is a basis of multilinear polynomials for the $T_{G/\langle g \rangle}$-ideal $Id_{G/\langle g \rangle}(A)$,  with the grading induced by $q$, then the set $$S=\{f\in K\langle X_G \rangle \hspace{0,1cm}|\hspace{0,1cm} Q(f)\in S^{\prime} \}$$ is a basis for the $T_G$-ideal $Id_{G}(A)$.
\end{proposition}
\textit{Proof.}

Let $f(x_{1,g_1},\dots, x_{n,g_n})$ be a multilinear polynomial. Denote by $a$ the invertible element in $A_g$. A homogeneous element of degree $q(g)$ may be written as $ba^{k}$ for some $b\in A_g$ and some integer $k$. Therefore any $Q(f)$-admissible substitution is of the form $(b_1a^{k_1},\cdots, b_na^{k_n})$, where $k_1,\dots, k_n$ are integers and $(b_1,\dots, b_n)$ is an $f$-admissible substitution. Since $a$ is central it is clear that $Q(f)(b_1a^{k_1},\cdots, b_na^{k_n})=a^kf(b_1,\dots, b_n)$. Hence $Q(f)$  is a graded identity for $A$ with the grading induced by $q$ if and only if $f$ is a $G$-graded identity for $A$. Now assume that $f\in Id_G(A)$, we may write $$Q(f)=\sum_{i=1}^n g_i s_i^{\prime}(w_i^{1},\dots, w_i^{k_i})h_i,$$ where $s_i^{\prime}\in S^{\prime}$, $g_1,h_i, w_i^{1},\dots,w_i^{k_i}\in K\langle X_{G/H} \rangle$. Now let $R:K\langle X_{G/H} \rangle \rightarrow K\langle X_{G} \rangle$ be a homomorphism such that $R(x_{i,q(g_i)})=x_{i,g_i}$ for $i=1,\dots, n$. Then it is clear that $R(Q(f))=f$ and therefore we have $$f=\sum_{i=1}^n R(g_i) R(s_i^{\prime})(R(w_i^{1}),\dots, R(w_i^{k_i}))R(h_i),$$ and since $Q(R(s_i))=s_i$ it follows that $R(s_i)\in S$. We have proved that any multilinear $G$-graded identity for $A$ lies in  the $T_G$-ideal  $\langle S \rangle_{T_G}$, generated by $S$. Hence $Id_G(A)\subseteq \langle S \rangle_{T_G}$. It is clear that $S\subseteq Id_G(A)$, this implies the inclusion $\langle S \rangle_{T_G}\subseteq Id_G(A)$. 
\hfill $\Box$

\begin{remark}
Note that if we assume, in the previous proposition, that $S^{\prime}$ is a basis of multilinear polynomials for the $T_{G/\langle g \rangle}$-space $C_{G/\langle g \rangle}(A)$, with the induced grading, then a simple modification of the proof yields that $S$ is a basis for the $T_G$-space $C_G(A)$.
\end{remark}

Let $G=(a)_4$ and $\Gamma$ the grading by the group $G$ on $M_2(\mathbb{C})$ that corresponds to $M_2(\mathbb{C}, \mathbb{Z}_4)$. The algebra $M_2(\mathbb{C}, \mathbb{Z}_4)$ has an invertible central element of degree $a^{2}$. The homogeneous components of grading on $S=M_2(\mathbb{C})$ induced by $q$ are the complex subspaces generated by $I$ and $C$ (neutral component) and the complex subspace generated by $A$ and $B$. It is clear that this is isomorphic as a graded algebra to $M_2^{(2)}\otimes \mathbb{C}$ and therefore satisfies the same graded identities as $M_2^{(2)}$. Hence the previous results yield the following corollary.

\begin{corollary}
Let $G=(a)_4$ and let $R$ be the algebra $M_2(\mathbb{C})$ with the $G$-grading weakly isomorphic to $M_2(\mathbb{C}, \mathbb{Z}_4)$. 
\begin{enumerate}
\item[(i)] A basis for $T_G(R)$ consists of the polynomials \[x_{1 k}x_{2 k}-x_{2 k}x_{1 k} and x_{1h}x_{2h}x_{3h}-x_{3h}x_{2h}x_{1h},\] where $k\in \{e, a^2\}$ and $h\in \{a, a^3\}$.
\item[(ii)] A basis for $C_G(R)$ consists of the polynomials \[x_{1h}^2x_{2h}+x_{2h}^2x_{1h}, \hspace{0,3cm} x_{1g}[x_{2k},  x_{3k}]x_{4g}, \hspace{0,3cm} x_{1g}\left(x_{2h}x_{3h}x_{4h}-x_{4h}x_{3h}x_{2h}\right)x_{5g},\] where $k\in \{e, a^2\}$, $h\in \{a, a^3\}$ and $g\in G$.
\end{enumerate}
\end{corollary}

\subsection{Non-regular Pauli Gradings}\label{NRPG}
Let $A$ denote the real algebra $M_n(\mathbb{C})$ with a non-regular Pauli grading by the finite abelian group $G$ and let $\beta_A:G\times G \rightarrow \mathbb{C}$ be the associated complex bicharacter. Complex division gradings on $M_n(\mathbb{C})$ have been classified and the bicharacter $\beta_A$ is non-degenerate (see \cite{BZ}).

\begin{remark}
Note that $f$ is a proper central polynomial for $A$ if and only if $f$ is not a graded identity and $\mathrm{deg}_G\ (f)=e$. Hence if $S$ is a basis for $Id_G(A)$ then $\{x_{1e}\}\cup S^{\prime}$ is a basis for $C_G(A)$, here $S^{\prime}$ is the set of polynomials $x_{ig}fx_{jh}$ where $f\in S$, $g,h \in G$ and $i,j$ are such that $x_{ig}$ and $x_{jh}$ do not appear in $f$.
\end{remark}

Given $\textbf{g}=(g_1,\dots,g_n)$ an $n$-tuple of elements of $G$ and a permutation $\sigma$ in $S_n$ we recall that $\gamma_{\sigma}^{\textbf{g}}$ is the  complex number such that \[a_1\cdots a_n= \gamma_{\sigma}^{\textbf{g}}\cdot a_{\sigma(1)}\cdots a_{\sigma(n)},\] for every $a_1\in A_{g_1},\dots, a_{n}\in A_{g_n}$. If $\gamma_{\sigma}^{\textbf{g}}$ is a real number then
\[
x_{i_1,g_{1}}\cdots x_{i_n,g_{n}}-\gamma_{\sigma}^{\textbf{g}}x_{i_{\sigma(1)},g_{\sigma(1)}}\cdots x_{i_{\sigma(n)},g_{\sigma(n)}}, \tag{I} \label{idrealgen}
\]
is a graded identity for $A$ for any $n$-tuple $(i_1,\dots, i_n)$ of natural numbers. If $\sigma, \tau$ are permutations such that $\gamma_{\sigma}^{\textbf{g}},\gamma_{\tau}^{\textbf{g}}\in \mathbb{C}\setminus \mathbb{R}$ are linearly independent over $\mathbb{R}$ then let $p,q \in \mathbb{R}$ such that $p(\gamma_{\sigma}^{\textbf{g}})^{-1}+q(\gamma_{\tau}^{\textbf{g}})^{-1}+1=0$. The polynomial
\[
x_{i_1,g_{1}}\cdots x_{i_n,g_{n}}+px_{i_{\sigma(1)},g_{\sigma(1)}}\cdots x_{i_{\sigma(n)},g_{\sigma(n)}}+qx_{i_{\tau(1)},g_{\tau(1)}}\cdots x_{i_{\tau(n)},g_{\tau(n)}}, \tag{II} \label{idcompgen}
\]
is a graded identity for $A$ for any $n$-tuple $(i_1,\dots, i_n)$ of natural numbers.

\begin{proposition}\label{lincomb}
If a polynomial is graded identity for $A$ then it is a linear combination of the polynomial identities in (\ref{idrealgen}) and (\ref{idcompgen}).
\end{proposition}
\textit{Proof.}
Let $f$ be a graded identity for $A$. We may assume that $f=\sum_{\sigma \in S_n} \mu_{\sigma}m_{\sigma}$, where $\mu_{\sigma}$ are real scalars and $$m_{\sigma}=x_{i_{\sigma(1)},g_{\sigma(1)}}\cdots x_{i_{\sigma(n)},g_{\sigma(n)}}.$$ We denote $\textbf{g}=(g_{1},\dots, g_{n})$ and assume that $\mu_{id}=1$. Denote $s$ the number of non-zero scalars $\mu_{\sigma}$. The result is proved by induction on $s$. Since $A$ has no monomial identities we have $s>1$. If $s=2$ then $f$ is of the form (\ref{idrealgen}). If $\gamma_{\sigma}^{\textbf{g}}\in \mathbb{R}$ for some $\sigma \in S_n$ such that $\lambda_{\sigma}\neq 0$ then $$f-\left(x_{i_1,g_{1}}\cdots x_{i_n,g_{n}}-\gamma_{\sigma}^{\textbf{g}}x_{i_{\sigma(1)},g_{\sigma(1)}}\cdots x_{i_{\sigma(n)},g_{\sigma(n)}}\right)$$ is an identity with $s-1$ or $s-2$ non-zero scalars, hence the result follows from the induction hypothesis. If $\gamma_{\sigma}^{\textbf{g}}\in \mathbb{C}\setminus \mathbb{R}$ for every $\sigma \in S_n$ such that $\lambda_{\sigma}\neq 0$ then $m>2$. Let $\sigma,\tau \in S_n$ such that $\mu_{\sigma}$ and $\mu_{\tau}$ are not zero. If $\gamma_{\sigma}^{\textbf{g}}$ and $\gamma_{\tau}^{\textbf{g}}$ are linearly dependent over $\mathbb{R}$ then there exists $p\in \mathbb{R}$ such that $m_{\sigma}-pm_{\tau}$ is a graded identity for $A$, this polynomial is as in (\ref{idrealgen}). In this case we apply the induction hypothesis to $$f-\mu_{\sigma}\left(m_{\sigma}-pm_{\tau}\right)$$ and the result follows. Otherwise there exists $p,q \in \mathbb{R}$ such that $p(\gamma_{\sigma}^{\textbf{g}})^{-1}+q(\gamma_{\tau}^{\textbf{g}})^{-1}+1=0$. The induction hypothesis applies to $$f-\left(x_{i_1,g_{1}}\cdots x_{i_n,g_{n}}+px_{i_{\sigma(1)},g_{\sigma(1)}}\cdots x_{i_{\sigma(n)},g_{\sigma(n)}}+qx_{i_{\tau(1)},g_{\tau(1)}}\cdots x_{i_{\tau(n)},g_{\tau(n)}}\right)$$ and the result follows.
\hfill $\Box$

The $T_G$-ideal $Id_G(A)$ admits a finite basis of identities (see \cite{AK}, \cite{S}), therefore there exists a finite basis for $Id_G(A)$ of polynomials in  (\ref{idrealgen}) and (\ref{idcompgen}). Our goal now is to exhibit such a finite set of identities.

Let $g, h \in G$, if $\beta_A(g,h)\in \mathbb{R}$ then
\begin{equation}\label{idreal}
x_{1,g}x_{2,h}-\beta_A(g,h)x_{2,h}x_{1,g},
\end{equation}
is a graded identity for $A$. If $\beta_A(g,h)\in \mathbb{C}\setminus \mathbb{R}$ then let $p,q \in \mathbb{R}$ be such that $\beta_A(g,h)$ is a root of $z^{2}+pz+q$. In this case
\begin{equation}\label{idcomp}
x_{1,g}x_{2,g}x_{3,h}+p(x_{1,g}x_{3,h}x_{2,g})+q(x_{3,h}x_{1,g}x_{2,g}),
\end{equation}
is a graded identity for $A$.
These are just particular polynomials of types (\ref{idrealgen}) and (\ref{idcompgen}).

\begin{lemma}\label{dist}
Let $A$ be an algebra with a non-regular Pauli grading such that for any $g,h \in G$ we have $\beta_A(g,h)\neq i$  and let $f(x_{1,g_1},\dots,x_{n,g_n})$ be a multilinear polynomial. If $f$ is  a graded identity for $A$ then there exists a multilinear polynomial $f^{\prime}(x_{1,g^{\prime}_1},\dots,x_{m,g^{\prime}_m})$ in $Id_G(A)$, where $g^{\prime}_1,\dots, g^{\prime}_m$ are pairwise distinct elements of $G$, such that $f$ lies in the $T_G$-ideal generated by $f^{\prime}$ together with the polynomials in (\ref{idreal}) and (\ref{idcomp}).
\end{lemma}

\textit{Proof.}
We prove the lemma by induction on the number $n$ of indeterminates, the case $n=1$ is obvious. If $g_1,\dots, g_n$ are pairwise distinct there is nothing to prove. Otherwise there exists indexes $i<j$ such that $g_i=g_j$, we may assume that $i=n-1$ and $j=n$. Let $m$ be a monomial in $f$, we claim that there exists a multilinear polynomial $f_{m}(x_{1,g_1},\dots,x_{n-2,g_{n-2}},x_{n-1,g^{\prime}})$, where $g^{\prime}=g_n^2$, such that $$m-f_{m}(x_{1,g_1},\dots,x_{n-2,g_{n-2}},(x_{n-1,g_{n-1}}\cdot x_{n,g_{n}}))$$ lies in the $T_G$-ideal generated by (\ref{idreal}) and (\ref{idcomp}). Let $w_1$, $w_2$ and $w_3$ be monomials such that
\begin{equation}
m=w_1x_{n-1,g}w_2x_{n,g}w_3,
\end{equation}
where $g=g_{n-1}=g_n$.
Denote $h$ the degree of the monomial $w_2$. If $\beta_A(g,h)\in \mathbb{R}$ then $$w_1x_{n-1,g}w_2x_{n,g}w_3-\beta_A(g,h)w_1w_2x_{n-1,g}x_{n,g}w_3$$ is a consequence of a polynomial in (\ref{idreal}), in this case $f_m=\beta_A(g,h)w_1w_2x_{n-1,g^{\prime}}w_3$. If $\beta_A(g,h)\in \mathbb{C}\setminus \mathbb{R}$ let $p,q$ be the real numbers such that $\beta_A(g,h)$ is a root of $z^{2}+pz+q$. Since $\beta_A(g,h)$ is a root of the unit different from $i$ then it follows that $p\neq 0$. In this case let $$f_m=\lambda_1w_1x_{n-1,g^{\prime}}w_2w_3+\lambda_2w_1w_2x_{n-1,g^{\prime}}w_3,$$ where $\gamma_1=-\frac{1}{p}$, $\gamma_2=-\frac{q}{p}$. It is clear that $$m-f_{m}(x_{1,g_1},\dots,x_{n-2,g_{n-2}},(x_{n-1,g}\cdot x_{n,g}))$$ is a consequence of (\ref{idcomp}). This proves the claim. Now write $$f=\sum_i\lambda_i m_i,$$ where $m_i$ are monomials and $\widehat{f}=\sum_i\lambda_if_{m_i}$. The difference $$f(x_{1,g_1},\dots,x_{n-2,g_{n-2}},x_{n-1,g_{n-1}}, x_{n,g_{n}})-\widehat{f}(x_{1,g_1},\dots,x_{n-2,g_{n-2}},(x_{n-1,g_{n-1}}\cdot x_{n,g_{n}}))$$ lies in the $T_G$-ideal generated by  (\ref{idreal}) and (\ref{idcomp}). We conclude that $$\widehat{f}(x_{1,g_1},\dots,x_{n-2,g_{n-2}},(x_{n-1,g_{n-1}}\cdot x_{n,g_{n}}))$$ is a graded identity for $A$. Since $A$ has a Pauli grading we conclude that the multilinear polynomial $\widehat{f}$ is a graded identity for $A$. Our induction hypothesis implies that there exists a graded identity $f^{\prime}(x_{1,g^{\prime}_1},\dots,x_{m,g^{\prime}_m})$, where $g^{\prime}_1,\dots, g^{\prime}_m$ are pairwise distinct, such that $\widehat{f}$ lies in the $T_G$-ideal generated by $f^{\prime}$ together with the identities in (\ref{idreal}) and (\ref{idcomp}). Clearly $f$ also lies in this $T_G$-ideal.
\hfill $\Box$

\begin{theorem}\label{noti}
Let $A$ be an algebra with a non-regular Pauli grading such that $\beta_A(g,h)\neq i$ for all $g,h \in G$. The set $S$ of the polynomials in (\ref{idreal}) and (\ref{idcomp}), together with the polynomials in (\ref{idrealgen}) and (\ref{idcompgen}) such that $g_1,\dots, g_n$ are pairwise distinct elements of $G$, is a basis for $Id_G(A)$.
\end{theorem}
\textit{Proof.}
Denote $I$ the $T_G$-ideal generated by $S$. Since $S\subseteq Id_G(A)$ it follows that $I\subseteq Id_G(A)$. The reverse inclusion follows from Lemma \ref{dist} and Propositon \ref{lincomb}.
\hfill $\Box$

We now consider the case where $\beta_A(g,h)=i$ for some $g,h \in G$. It is clear that in this case the polynomial
\begin{equation}\label{reorder}
x_{1,g}x_{2,h}x_{3,g}-x_{3,g}x_{2,h}x_{1,g},
\end{equation}
is a graded identity for $A$.

\begin{remark}
The polynomial in (\ref{reorder}) is a graded identity for $A$ even if $\beta_A(g,h)\neq i$. In this case it is a consequence of the identities (\ref{idreal}) and (\ref{idcomp}).
\end{remark}

 Let $g,h_1,h_2,h_3$ be elements of $G$ such that $\beta_A(g,h_i)=i$ for $i=1,2,3$. The polynomial
\begin{equation}\label{idi}
x_{1,g}x_{2,h_1}x_{3,g}x_{4,h_2}x_{5,g}x_{6,h_3}x_{7,g}+x_{1,g}x_{3,g}x_{5,g}x_{7,g}x_{2,h_1}x_{4,h_2}x_{6,h_3}
\end{equation}
is a graded identity for $A$.

\begin{lemma}\label{lim}
Let $A$ be an algebra with a non-regular Pauli grading by the group $G$ and let $f(x_{1,g_1},\dots,x_{n,g_n})$ be a multilinear polynomial. If $f$ is  a graded identity for $A$ then there exists a multilinear polynomial $f^{\prime}(x_{1,g^{\prime}_1},\dots,x_{m,g^{\prime}_m})$ in $Id_G(A)$, where for each $g\in G$ the number of indexes $j$ such that $g_j^{\prime}=g$ is at most $3$, such that $f$ lies in the $T_G$-ideal generated by $f^{\prime}$ together with the polynomials in (\ref{idreal}), (\ref{idcomp}), (\ref{idi}) and (\ref{reorder}).
\end{lemma}

\textit{Proof.}
The proof is analogous to the proof of Lemma \ref{dist}. If $g_1,\dots, g_n$ are such that for each $g\in G$ the number of indexes $i$ such that $g_i=g$ is at most $3$ there is nothing to prove. Otherwise there exists indexes $i<j<k<l$ such that $g_i=g_j=g_k=g_l$, we denote $g$ this element. Let $m$ be a monomial in $f$ and $w_i$, $i=0,\dots 4$, monomials such that
\begin{equation}
m=w_0x_{i,g}w_1x_{j,g}w_2x_{k,g}w_3x_{l,g}w_4.
\end{equation}

Denote $h_i$ the degree of the monomial $w_i$. We may assume without loss of generality that $\{n,n-1,n-2,n-3\} =\{i,j,k,l\}$. If $\beta_A(g,h_1)\neq i$ we use the identities (\ref{reorder}) to assume that $i=n-1$ and $j=n$, then as in the proof of Lemma \ref{dist} determine a polynomial $f_{m}(x_{1,g_1},\dots,x_{n-2,g_{n-2}},x_{n-1,g^2})$, such that
\begin{equation}\label{fm}
m-f_{m}(x_{1,g_1},\dots,x_{n-2,g_{n-2}},(x_{n-1,g_{n-1}}\cdot x_{n,g_{n}}))
\end{equation}
lies in the $T_G$-ideal generated by (\ref{idreal}) and (\ref{idcomp}).  Analogously if $\beta_A(g,h_j)\neq i$ for $j=2$ or $j=3$ we obtain a multilinear polynomial $f_{m}(x_{1,g_1},\dots,x_{n-2,g_{n-2}},x_{n-1,g^2})$ such that $m$ lies in the $T_G$-ideal  generated by $f_m$ and the polynomials in (\ref{reorder}), (\ref{idreal}) and (\ref{idcomp}). Finally we assume that $\beta_A(g,h_j)=i$ for $j=1,2,3$, in this case let
$$f_{m}(x_{1,g_1},\dots,x_{n-2,g_{n-2}},x_{n-1,g^2})=-w_0x_{n-3,g}x_{n-2,g}x_{n-1,g^2}w_1w_2w_3w_4.$$ It is clear that, for this polynomial, the difference in (\ref{fm}) lies in the $T_G$-ideal generated by the polynomials in (\ref{reorder}) and (\ref{idi}).

Now write $$f=\sum_i\lambda_i m_i,$$ where $m_i$ are monomials and $\widehat{f}=\sum_i\lambda_if_{m_i}$. The multilinear polynomial  $\widehat{f}$ is a graded identity for $A$ in $n-1$ indeterminates and the result follows by induction on $n$.
\hfill $\Box$

The proof of Theorem \ref{noti} yields mutatis mutandis (with the previous lemma instead of Lemma \ref{dist}) the following result:

\begin{theorem}
Let $A$ be an algebra with a non-regular Pauli grading. The set $S$ of the polynomials in (\ref{idreal}), (\ref{idcomp}) (\ref{idi}) and (\ref{reorder}), together with the polynomials in (\ref{idrealgen}) and (\ref{idcompgen}) such that for any $g\in G$ the number of indexes $j$ such that $g_j=g$ is at most $3$, is a basis for $Id_G(A)$.
\end{theorem}

\section{Basis for the Graded Identities and Central Polynomials for Real Graded Division Algebras}\label{gdiv}

In this section we consider the real graded division algebras classified in \cite{BZ2}. The classification is in terms of four basic graded division algebras that we describe next.

\subsection{Basic Real Graded Division Algebras}\label{basic2}

Let $G$ be a finite abelian group, a $2$-cocycle is a map $\sigma:G\times G \rightarrow K^{\times}$ such that $$\sigma(g,h) \sigma(gh,k)=\sigma(g,hk)\sigma(g,k).$$ In this case $\beta(g,h):=\sigma(g,h)\sigma(h,g)^{-1}$ is a skew-symmetric bicharacter, conversely every skew-symmetric bicharacter is obtained from a $2$-cocycle in this way (see \cite{Sch}, \cite[Remark 13]{AO}). Denote $\mathbb{C}^{\sigma}G$ the (complex) twisted group algebra, that is, the complex vector space with basis $G$ and multiplication such that $g\cdot h:= \sigma(g,h)gh$. This algebra with its canonical $G$-grading is denoted $P(\beta)$, moreover this algebra viewed as an algebra over $\mathbb{R}$ is denoted $P(\beta)_{\mathbb{R}}$. 

Now let $m>1$ be a natural number and $\epsilon \in \{1,-1\}$, we denote $\mathcal{D}(m,\epsilon)$ the subalgebra in the complex group algebra $\mathbb{C}(g)_m$ (of the cyclic group of order $m$) generated by $\mu g$, where $\mu^{m}=\epsilon$. This is a commutative $(g)_m$-graded division algebra. We remark that $\mathcal{D}(2,-1)$ is weakly isomorphic to $\mathbb{C}^{(2)}$.

Let $G=(g)_k\times (h)_{l}$ be the product of two cyclic $2$-groups. Let $\mu, \nu \in \{1,-1\}$, we fix complex numbers $\epsilon$, $\eta$ such that $\epsilon^k=\mu$, $\eta^l=\nu$ and denote $\mathrm{D}(k,l;\mu,\nu)$ the subalgebra of $S=(\mathbb{R}G \otimes M_2(\mathbb{C}))$ generated by $u=g\otimes \epsilon A$ and $v=h\otimes \eta B$. We recall that $
A=\left(\begin{array}{cc}
	1 & 0\\
	0 & -1
\end{array}\right)
$ and $
B=\left(\begin{array}{cc}
	1 & 0\\
	0 & -1
\end{array}\right)
$ are the Pauli matrices introduced in Section \ref{divgrad}. The algebra $S$ has a canonical grading where $S_g=(\mathbb{R}g)\otimes M_2(\mathbb{C})$, it is clear that $\mathcal{D}(k,l;\mu,\nu)$ is a homogeneous subalgebra. 

Let $n$ be a natural number and $\epsilon\in \{1,-1\}$. We denote $\mathcal{E}(\epsilon,n)$ the subalgebra of $\mathbb{R}(g)_n\otimes M_2(\mathbb{C})$ generated by $u=1\otimes C$ and $v=g\otimes \epsilon A$. We recall that $C=AB$.

\begin{theorem}\cite[Theorem 9.1]{BZ2}
Let $\Gamma: R=\oplus_{g\in G}$ be a grading of finite-dimensional
real graded division algebra $R$, $G = Supp \Gamma$. Then one of the following are true.
\begin{enumerate}
\item[(i)] If $R$ is commutative with $\mathrm{dim}\ R_e=1$ then $\Gamma$ is equivalent to the graded
tensor product of algebras $\mathcal{D}(m;\eta)$.
\item[(ii)] If $\mathrm{dim}\ R_e=1$ and $R$ is not commutative, then $\Gamma$ is equivalent to the
graded tensor product of several copies of $\mathcal{D}(2^k, 2^l; \mu, \nu)$ and several copies
of $\mathcal{D}(m;\eta)$, where $m$, $k$, $l$ are natural numbers, $m>1$ and $\mu, \nu, \eta = \pm1$.
\item[(iii)] If $\mathrm{dim}\ R_e=4$ then $\Gamma$ is isomorphic to the graded tensor product of the
trivial grading on $\mathbb{H}$, several copies of $\mathcal{D}(2^k, 2^l; \mu, \nu)$ and several copies of
$\mathcal{D}(m;\eta)$, where $k$, $l$, $m$ are natural numbers $m > 1$ and $\mu, \nu, \eta = \pm1$.
\item[(iv)] If $\mathrm{dim}\ R_e=2$ and $R_e$ is not central then $\Gamma$ is equivalent to the tensor
product of one copy of $\mathcal{E}(\epsilon,n)$, several copies of $\mathcal{D}(2^k, 2^l; \mu, \nu)$ and several
gradings of the form $\mathcal{D}(m;\eta)$ where $m$, $k$, $l$ are natural numbers $m > 1$, $n$
a 2-power, and $\epsilon, \mu, \nu, \eta = \pm1$.
\item[(v)] If $\mathrm{dim}\ R_e=2$ and $R_e$ is central then $\Gamma$ is isomorphic to $P(\beta)_{\mathbb{R}}$ for a
skew-symmetric complex bicharacter $\beta:G\times G \rightarrow \mathbb{C}^{\times}$.
\end{enumerate}
\end{theorem}

\subsection{Basis For the Graded Identities and Central Polynomials}\label{basis2}

The algebras in $(i)$ and $(ii)$ have regular gradings, basis for the graded identities and central polynomials are obtained as in Section \ref{reg}, the algebras in $(iii)$ have a tensor product grading of $\mathbb{H}$ with the trivial grading and an algebra with a regular grading, in this case basis are obtained as in Section \ref{TPDG}, the algebras in $(v)$ that are not regular will have basis described in Section \ref{NRPG}. To determine basis for the algebras in $(iv)$ we need only to consider the algebras $\mathcal{E}(\epsilon, 2^k)$. 

\begin{proposition}
Let $G=(g)_n$ be the cyclic group of order $n$, where $n=2^k$ for some $k$, $H=\langle g^2 \rangle$ and $q:G\rightarrow G/H$ the canonical homomorphism. If a polynomial $f$ in $K\langle X_{G/H} \rangle$ is multilinear then $f$ is a graded identity (resp. graded central polynomial) for $\mathcal{E}(\epsilon, 2^k)$, with the grading induced by $q$, if and only if $f$ is a graded identity (resp. graded central polynomial)for $M_2(\mathbb{R})$ with the $G/H$-grading weakly isomorphic to $M_2^{(2)}$. 
\end{proposition}
\textit{Proof.}
The element $w=g\otimes \epsilon I$ is an invertible central element. Let $R=\mathcal{E}(\epsilon, 2^k)$ with the grading induced by $q$ and denote $G^{\prime}=\{e, a\}$ the group $G/H$. The set $$\beta=\{w^{2m}(1\otimes I), w^{2m}(1\otimes C), w^{2m-1}(1\otimes A), w^{2m-1}(1\otimes B)| m=0,1,\dots,2^{k-1}\}$$ is a basis for $R$ of homogeneous elements. Let $f(x_{1g_1^{\prime}},\dots, x_{ng_n^{\prime}})$ be a multilinear polynomial in $K\langle X_{G/H} \rangle$, this polynomial is a graded identity for $R$ if and only if the result of every admissible substitution by elements in $\beta$ is zero. Let $(a_1,\dots,a_n)$ be an admissible substitution by elements of $\beta$, then $a_i=w^{p_i}(1\otimes b_i)$, where $b_i\in \{I, A, B, C\}$ and $p_i$ is a natural number. In this case $(b_1,\dots, b_n)$ is an admissible substitution by elements of $M_2(\mathbb{R})$ with the $G/H$-grading weakly isomorphic to $M_2^{(2)}$. We have $f(a_1,\dots,a_n)=w^{p}(1\otimes f(b_1,\dots, b_n))$, where $p=p_1+\cdots + p_n$. Since $w$ is an invertible central element we conclude that $f$ is a graded identity (resp. graded central polynomial) for $R$ if and only if it is a graded identity (resp. graded central polynomial) for $M_2(\mathbb{R})$.
\hfill $\Box$

Let $v=g\otimes \epsilon A$, the element $v^2=-g^2\otimes \epsilon^2I$ is a homogeneous invertible central element of degree $g^2$. Therefore the previous proposition together with Proposition \ref{1}, Lemma \ref{2} and Theorem \ref{3} proves the following:

\begin{corollary}
Let $G=(g)_{2^k}$ and $R$ the $G$-graded division algebra $\mathcal{E}(\epsilon, 2^k)$. 
\begin{enumerate}
\item[(i)] A basis for $T_G(R)$ consists of the polynomials \[x_{1 k}x_{2 k}-x_{2 k}x_{1 k}\mbox{ \hspace{0,3cm } and \hspace{0,3cm }}x_{1h}x_{2h}x_{3h}-x_{3h}x_{2h}x_{1h},\] where $k\in \{g^{2m}| m=0,1,\dots, 2^{k-1}\}$ and $h\in \{g^{2m+1}| m=0,1,\dots, 2^{k-1}\}$.
\item[(ii)] A basis for $C_G(R)$ consists of the polynomials \[x_{1h}^2x_{2h}+x_{2h}^2x_{1h}, \hspace{0,3cm} x_{1g^{\prime}}[x_{2k},  x_{3k}]x_{4g^{\prime}}\mbox{ \hspace{0,3cm } and \hspace{0,3cm }}x_{1g^{\prime}}\left(x_{2h}x_{3h}x_{4h}-x_{4h}x_{3h}x_{2h}\right)x_{5g^{\prime}},\] where $k\in \{g^{2m}| m=0,1,\dots, 2^{k-1}\}$, $h\in \{g^{2m+1}| m=0,1,\dots, 2^{k-1}\}$ and $g^{\prime}\in G$.
\end{enumerate}
\end{corollary}

\end{document}